\documentclass[12pt]{amsart}
\usepackage{german}
\usepackage[ansinew]{inputenc}
\usepackage{hyperref}\usepackage{graphicx}

\newcommand{\tang}{\operatorname{tang}}
\title{Ein Brief von Eisensteins Eltern an Gauss}
\author{Franz Lemmermeyer}
\begin{document}
\maketitle

Constantin und Helene Eisenstein, die Eltern von Gotthold Eisenstein,
haben sich nach dem Tod ihres Sohnes brieflich bei Gauss bedankt.
Dieser Brief ist inzwischen online\footnote{Siehe
  \url{https://gauss.adw-goe.de/handle/gauss/2588?}.}
gestellt. Hier wollen wir eine Transkription des Briefes\footnote{Erstellt
  mit freundlicher Unterstützung durch Menso Folkerts, dem ich dafür
  herzlich danke. Ebenso herzlich danke ich Peter Ullrich für Kommentare
  und Korrekturen.} geben, sowie den mathematisch interessanten
Teil der Anlage, nämlich einen bisher unveröffentlichten Beweis des
quadratischen Reziprozitätsgesetzes mit Hilfe der Tangensfunktion.

\section{Der Brief von Eisensteins Eltern an Gauß}
  
An \\
S[eine]r Hochwohlgeboren \\
den Herrn Geheimen-Hofrath \\
Professor Dr Gauss \\
Zu Göttingen \\

\hfill   	Berlin den 12tn November 1852.

Hochwohlgeborner Herr!
\medskip

Wenn gleich Herr Professor D$^{\text{r}}$ Encke die Güte gehabt hat
Ihnen Hochgeehrter
Herr das Ableben unseres Sohnes Dr Gotthold Eisenstein zur Zeit anzuzeigen;
so bewegt uns doch die Pflicht der Dankbarkeit, einige Zeilen an den
Wohlthäter und Meister unseres Kindes Ehrfurchtsvoll zu richten.

Mit traurigem Herzen wiederholen wir Ihnen, daß unser vielgeliebtes
letztes Kind (dem Fünf voraus gegangen sind) Ferdinand Gotthold Maximilian
Eisenstein Dr der Phylosophie etc geboren zu Berlin am 16tn April 1823,
am 11tn October 1852. früh 1/4 vor 6 Uhr in Folge von Nervenabzehrung
gestorben ist. Unser Schmerz ist groß! Doch hat der Allmächtige wohlgethan,
denn er war sein ganzes Leben hindurch leidend und war es uns Eltern von
der Vorsehung nicht beschieden ihm ein sorgenfreies Leben bereiten zu
können, wodurch er auch von dieser Seite stets bedrängt war.

Die väterliche Theilnahme die Sie Hochgeehrter Meister an unserm Sohne
persönlich und die Anerkenntniße an seinen Leistungen haben ihn oft in
trüben Stunden aufgerichtet und wieder frisch an's Werk gehen
lassen. Wir statten Ihnen dafür unsern herzlichsten Dank ab.

Welcher Idial Meister Sie ihm waren, läßt sich mit Worten nicht
aussprechen! Ihre liebevollen Briefe die Sie an ihn gerichtet, sind in
unsern Händen und sollen uns als theures Andenken von Ihnen
Hochgeehrter Herr verbleiben.

Das beikommende Buch hatte für ihn keinen anderen Namen als mein
\underline{Gauss}, daher bitten wir Sie, dasselbe mit seinen
Anmerkungen und Einlagen, wie wir es nach seinem Ableben vorgefunden
haben, hochgeneigst als Andenken an ihn annehmen zu wollen.
Auch erlauben wir uns, Ihnen seine letzte Academische Abhandlung und
Antrittsrede beizulegen, wie auch einige mathematische Arbeiten welche
seine letzten waren, und wenige Wochen vor seinem Tode ihm aus der Feder
geflossen sind. Wir haben in seinem Nachlasse sehr viele mathematische
Arbeiten gefunden in welchen manche schöne Gedanken enthalten sein dürften.

Möge der Allmächtige Ewr. Hochwohlgeboren noch lange in Ihren
Wirkungskreise gesund erhalten dieses wünschend von ganzem Herzen
\medskip

die Sie Hochverehrende

Eltern des Verstorbenen
\medskip

Constantin Eisenstein und

Helene Eisenstein

Neue Grünstraße N$^{o}$ 13

\section*{Der Beweis des quadratischen Reziprozitätsgesetzes}

Den analogen Beweis für das die quadratischen Reste in der reellen Theorie
betreffende Fundamentaltheorem mit Hülfe der Kreisfunktionen will ich hier
noch beifügen.
\begin{align*}
  (\cos x + i \sin x)^q & = \cos qx + i \sin qx
  = \sum_{m=0}^{m=q} i^m q_m \cos x\ {}^{q-m} \sin x\ {}^m \\
  \tang qx & = \frac{\sum_{n=0}^{n = \frac{q-1}2} (-1)^n q_{2n+1} \cos x\ {}^{q-2n-1}
             \sin x\ {}^{2n+1}}
                   {\sum_{n=0}^{n = \frac{q-1}2} (-1)^n q_{2n} \cos x\ {}^{q-2n}
                     \sin x\ {}^{2n}} \\
                   & =  \frac{\sum (-1)^n q_{2n+1} \tang x\ {}^{2n+1}}
                             {\sum (-1)^n q_{2n}   \tang x\ {}^{2n}}, 
  \end{align*}
welches sich, wenn $q$ eine ungerade Primzahl bedeutet, auf die Form
bringen läßt:
$$ \frac{q \phi(\tang x\ {}^2) + (-1)^{\frac{q-1}2} \tang x\ {}^{q-1}}
        {1 + q \psi(\tang x\ {}^2)} \tang x, $$
wo $\phi$ und $\psi$ ganze Funktionen von $\tang^2 x$ mit ganzen
Coefficienten sind. 

Sei $p$ eine von $q$ verschiedene ungerade Primzahl, dann hat man ebenso
$$ \tang px = \tang x \
  \frac{p - \frac{p(p-1)(p-2)}{1 \cdot 2 \cdot 3} \tang x\ {}^2
       + \ldots + (-1)^{\frac{p-1}2} \tang x\ {}^{p-1}}
    {1 -  \frac{p(p-1)}{1 \cdot 2} \tang x\ {}^2 + \ldots +
            (-1)^{\frac{p-1}2} p \tang x\ {}^{p-1}}. $$

Die Wurzeln der Gleichung
$$ (-1)^{\frac{p-1}2} Z^{p-1} + \ldots + p = 0 $$
sind dann offenbar in der Formel enthalten $\tang \frac{\rho \omega}{p}$,
wenn $\rho$ ein vollständiges Restsystem (mod $p$) mit Ausschluß der Null
durchläuft, und $\omega = 2\pi$ ist. Man hat hiernach
$$ \prod \tang \frac{\rho \omega}{p} = (-1)^{\frac{p-1}2} p. $$
(Die Multiplication links bezieht sich auf alle $\rho$.)

Die Werthe von $\rho$ lassen sich folgendermaßen gruppieren:
$$ \begin{array}{c|c}
  r_1 & -r_1 \\
  r_2 & -r_2 \\
  \vdots & \cdots \\
  r_{\frac{p-1}2} & - r_{\frac{p-1}2}
  \end{array} $$
Hieraus schließt man
\begin{equation}\label{Ea}
  \prod  \Big(\tang \frac{r \omega}p \Big)^2 = p. 
\end{equation}
$r$ bedeutet den Inbegriff der Zahlen $r_1$, $r_2$, \ldots, $r_{\frac{p-1}2}$.

Die Größen $(\tang \frac{r \omega}p )^2$ sind die Wurzeln einer Gleichung
mit ganzen Coefficienten
$$ (-1)^{\frac{p-1}2} Z^{p-1} + \ldots + p = 0, $$
also ist jede symmetrische Verbindung dieser Größen einer \underline{ganzen}
Zahl gleich.

Multipliciert man alle $r$ mit $q$, so hat man entweder
$qr \equiv r'$ oder $qr \equiv -r' \pmod p$, wo $r'$ sich
jedesmal unter den $r$ befindet. Nach diesen beiden Fällen hat man resp.
$$ \tang \frac{qr\omega}{p} = \tang \frac{r'\omega}p
   \quad \text{oder} \quad
   = - \tang \frac{r'\omega}p $$
also in beiden Fällen
$$ qr \equiv  r' \frac{\tang \frac{qr\omega}p}{\tang \frac{r'\omega}p}
\pmod p, $$
und da alle $r'$ alle $r$ erschöpfen, so schließt man hieraus
\begin{align*}
  q^{\frac{p-1}2} \prod (r) & \equiv \prod (r) \prod
  \bigg\{ \frac{\tang \frac{qr\omega}p}{\tang \frac{r\omega}p}\bigg\} \pmod p. \\
  \Big(\frac qp\Big) & =   \prod
  \bigg\{ \frac{\tang \frac{qr\omega}p}{\tang \frac{r\omega}p} \bigg\} \pmod p
\end{align*}
Aber nach dem Obigen ist
$$ \frac{\tang qx}{\tang x} =
   \frac{q \phi(\tang x\ {}^2) + (-1)^{\frac{q-1}2} \tang x\ {}^{q-1}}
        {1 + q \psi(\tang x\ {}^2)}, $$
also
\begin{align*}
  \prod \bigg\{ \frac{\tang \frac{qr\omega}p}{\tang \frac{r\omega}p} \bigg\}
  & = \prod  \frac{\big(q \phi(\{\tang \frac{r\omega}p\}^2)
        + (-1)^{\frac{q-1}2} (\tang \frac{r\omega}p)^{q-1}\big)}
        {1 + q \psi(\{\tang \frac{r\omega}p\}^2)} \\
        & = \frac{qP + (-1)^{\frac{p-1}2 \frac{q-1}2}
          \prod \{\tang \frac{r\omega}p\}^{q-1}}{1 + qQ} \\
        & = \frac{qP + (-1)^{\frac{p-1}2 \frac{q-1}2}
          p^{\frac{q-1}2}}{1 + qQ} \quad \text{nach (\ref{Ea})}
\end{align*}
also
$$ \Big(\frac qp \Big) =
  \frac{qP + (-1)^{\frac{p-1}2 \frac{q-1}2} p^{\frac{q-1}2}}{1 + qQ}, $$
  woraus unmittelbar das Reciprocitätsgesetz folgt.

  Die Werthe von $P$ und $Q$ können übrigens direct hingeschrieben werden,
  wenn man die symmetrischen Funktionen, welche sie darstellen, nach
  dem Newtonschen Theorem in den Coefficienten der Gleichung ausdrückt.
  Sie werden hiernach Aggregate von Binomialkoefficienten.

  Ähnliches gilt für das biquad. Fundamentaltheorem, nur kennt man hier 
  noch nicht das allgemeine Gesetz der Coefficienten in den
  Multiplicationsformeln. Letztere Coefficienten, die ich durch
  $A_1$, $A_2$, \ldots be\-zeich\-net habe, scheinen einiger Aufmerksamkeit
  würdig zu sein. Ich habe bloß gefunden, daß wenn man $m = a + bi$
  setzt, die in Rede stehenden Coefficienten durchaus nicht ganze Funktionen
  von $m$, sondern vielmehr Funktionen der beiden getrennten Variablen
  $a$ und $b$ sind.  Einiges über die Natur dieser Coefficienten
  kann man aus dem Umstande schließen, daß der $\mu$te Coeff. im Zähler
  und der $\mu$te im Nenner gleichzeitig verschwinden müssen für
  alle zusammengehörigen Werthe von $a$ und $b$, die $a^2 + b^2 < 4\mu + 1$
  machen.
  \medskip
  
  G. Eisenstein

  \medskip
  
  26. Januar, 45.

  Das cubische Reciprocitätsgesetz habe ich aus dem Integrale der
  Differentialgleichung
  $$ \frac{\partial y}{\sqrt{1-y^3}} =
  \Big\{a + b \frac{-1 + \sqrt{-3}}2 \Big\} \frac{\partial x}{\sqrt{1-x^3}} $$
  abgeleitet.

  \medskip
  
\includegraphics[width=11.4cm]{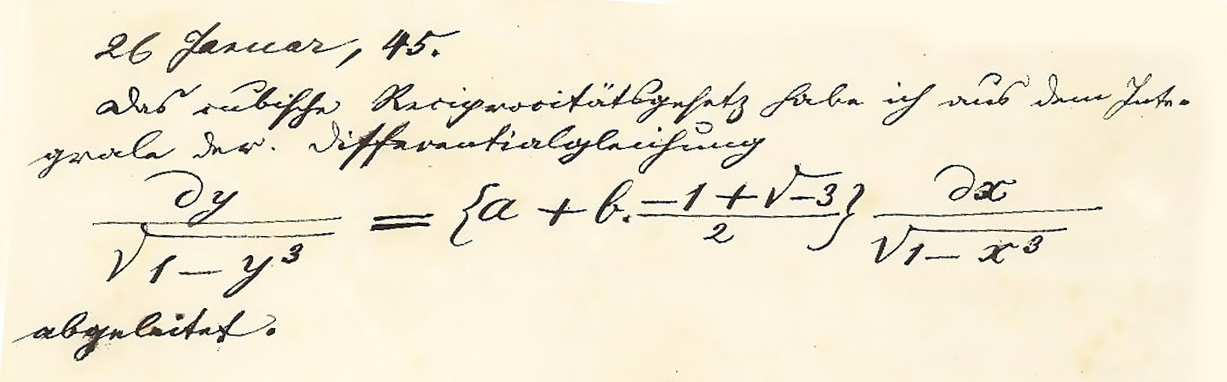}

\section{Kommentare}

Neben den oben transkribierten Teilen des Briefes samt Anlagen ist nur
noch der folgende kleine Ausschnitt interessant:

\begin{quote}
  {\em Ich darf fast gar nicht sprechen, weil ich zu sehr huste und seit
    8 Tagen schon 3 mal eine Menge Blut ausgeworfen habe. Sie können
    deshalb immer bei mir bleiben und mir etwas erzählen.

    Ich habe schon 1 Nacht draußen in Bethanien geschlafen, aber es war
    zu unruhig.  Mitten in der Nacht sprang ein Bierpropfen aus der
    Flasche etc Später jetzt bin ich zu schwach}
\end{quote}

Eisenstein wurde Ende Juli 1852, wie Biermann auf S.~927 von Eisensteins
Werken schreibt, nach einem Blutsturz für einige Tage in der Heilanstalt
Bethanien untergebracht. Vermutlich stammen diese Zeilen aus dieser Zeit.
An wen diese Zeilen adressiert waren, bleibt unklar.

\subsection*{Kommentare zum Brief}

Eisensteins Eltern haben Gauss als Andenken das Exemplar seiner Disquisitiones
geschenkt. Egon Ullrich, der 1925 in Graz bei Anton Rella promoviert hatte,
hat das Buch in der Bibliothek des mathematischen Instituts in Gießen
entdeckt; es gehörte zu den Büchern aus dem Nachlass von Eugen Netto.
1979 hat Benno Artmann diese Information an Andr\'e Weil weitergegeben,
der 1976 Eisensteins Werke besprochen und ein Buch über Eisensteins Zugang
zu elliptischen Funktionen veröffentlicht hatte. In \cite[S.~463]{Weil}
spekuliert er, dass Eisensteins Exemplar über den gemeinsamen Lehrer von
Eisenstein und Netto, Karl Heinrich Schellbach, an Netto gekommen sein könnte,
aber Peter Ullrich (siehe \cite[insbes. S.~206--207]{Ullr1} und \cite{Ullr2})
bemerkt, dass Netto das Buch von einem wenig mathematik-affinen Buchhändler
gekauft haben muss. Man wird wohl annehmen müssen, dass Ernst Schering, der
den Gaußschen Nachlass geordnet und die Herausgabe der Gaußschen Werke bis
zu seinem eigenen Tod geleitet hat, etwas damit zu tun gehabt hat. Warum
Schering, der zahlentheoretisch durchaus bewandert war und sich ausführlich
mit Beweisen des quadratischen Reziprozitätsgesetzes befasst hat, das
Eisensteinsche Exemplar weggegeben haben könnte, bleibt dabei ein Rätsel.
Auch was aus den mathematischen Arbeiten geworden ist, welche Eisensteins
Eltern dem Brief beigelegt haben, muss offen bleiben; die Antrittsrede zu
seiner Wahl in die Berliner Akademie, die Eisenstein am 1. Juli 1852
gehalten hat, ist in Band II seiner Werke abgedruckt.

Die angesprochenen Gaußschen \glqq{}Anerkenntniße an seinen Leistungen\grqq{}
waren für den jungen Eisenstein sicherlich eine große Motivation; auf der
anderen Seite haben diese bei seinen Zeitgenossen für manchmal unverhohlenen
Neid gesorgt. So schreibt Riemann am 23. Juli 1847 in einem Brief an
seinen Vater (\cite[S.~94]{Neuen}), Eisenstein habe sich \glqq{}bei Gauß
in Gunst zu setzen gewußt und er ist von diesem an Alexander von Humboldt
empfohlen.\grqq{}

Stein des Anstoßes war wohl der Brief vom 9. Juli 1845 von Gauß an von
Humboldt; dort erklärt Gauß, dass er zwar Dirichlet für den Orden
{\em Pour le m\'erite} vorgeschlagen habe, dass ihm diese Wahl aber
schwergefallen sei:
\begin{quote}
  {\em Sollte jene Injunction aber ganz buchstäblich genommen werden
    müs\-sen, nemlich unabhängig von jeder Rücksicht, also auch von
    der, ob einige Aussicht sei, daß dem Vorschlage, wie ungewöhnlich
    er auch sei, Folge gegeben werden könne, so bekenne ich, daß mir
    die Wahl zwischen Herrn Dirichlet und Herrn Eisenstein schwer geworden
    sein würde, da die Arbeiten des letztern in vollem Maße dasselbe
    Prädicament verdienen, wie die des erstern.}
\end{quote}
Jacobi lässt daraufhin seine alte Arbeit aus dem Jahre 1837 über Kreisteilung
in Crelles Journal wieder abdrucken und fügt eine Fußnote hinzu, in welcher
er Eisenstein des Plagiats bezichtigt: zum Einen soll er sich der
Kollegienhefte zu Jacobis Vorlesung über Zahlentheorie bedient haben, zum
andern stimme einer von Eisensteins Beweisen des quadratischen
Reziprozitätsgesetzes mit seinem eigenen überein, den Legendre in seine
Zahlentheorie aufgenommen habe. Auch Cauchy hat diesen Beweis (mehr oder
weniger zeitgleich mit Jacobi) veröffentlicht, und keinem der dreien ist
aufgefallen, dass es sich dabei im Wesentlichen um den sechsten Gaußschen
Beweis handelt; erst Gauß hat Eisenstein dies in seinem Brief vom
23. Juni 1844 an Encke (siehe \cite{Wittmann}) wissen lassen:
\begin{quote}
  {\em An der Art wie er sich über die bisherigen Beweise äußert,
    möchte ich fast vermuthen daß ihm das, was im XVII Band der
    Pariser Memoires S.~454 steht nicht bekannt geworden ist. Mir selbst
    ist dieser Band erst in diesen Tagen zu Gesicht gekommen, und die
    $3^{\text{te}}$ Edition von Legendre \glqq{}Th\'eorie des Nombres\grqq{},
    auf welche seine Stelle Bezug nimmt, habe ich überall noch nicht
    gesehen. Machen Sie doch Herrn Eisenstein auf diese Stelle aufmerksam,
  über welche ich übrigens seinem Urtheile nicht vorzugreifen brauche.}
\end{quote}

Offenbar hat es Jacobi nicht bei dieser Aktion belassen: In einem Brief
an Olfers\footnote{\cite[S. 105--106]{Olfers}.} schreibt Humboldt:
\begin{quote}
  {\em Da ich bestimmt weiß, daß Jacobi und selbst Encke auf die liebloseste
    Weise dem jungen blutarmen Mathematiker Eisenstein beim Minister
    Eichhorn zu schaden gesucht, so habe ich mir die süße Freude gemacht,
    beiden einen neuen Brief von Gauß zu schicken, in dem
    verbotenus\footnote{Lateinisch für wortwörtlich.} steht
    (14. April 1846):}
  \smallskip
  
      {\em Ob ich gleich voraussetzen darf, daß meine Empfehlung dem in Berlin
      lebenden jungen Mathematiker Doctor Eisenstein, eines Mannes, den
      ich sehr hochschätze, bei Ihrer Regierung ganz überflüssig ist,
      so will ich doch nicht unterlassen, es auszusprechen, daß ich seine
      Begabung wie eine solche betrachte, welche die Natur in jedem
      Jahrhundert nur wenigen erteilt.}
\end{quote}
Humboldt fährt fort, Jacobi habe 
\begin{quote}
  {\em seitdem ausgesprochen, daß Gauß ganz herunter und
  verstandesschwach geworden sei, daß Herr Eisenstein ja nicht daran
  denken müsse, sich hier als Privatdozent zu habilitieren, sondern
  nach Halle, oder nach Bonn gehen müsse!!!}
\end{quote}

\subsection*{Kommentar zum Beweis des Reziprozitätsgesetzes}

Der oben vorgestellte Beweis Eisensteins für das quadratische
Reziprozitätsgesetz unterscheidet sich nur formal von seinem bekannten
Beweis mit Hilfe der Sinusfunktion (siehe \cite{Eis1}). Der wesentliche
Unterschied der beiden Beweise liegt darin begründet, dass $\sin(qx)$
für ungerade Werte von $q$ ein Polynom in $\sin x$ ist, während $\tan(qx)$
eine rationale Funktion von $\tan x$ ist -- der Grund dafür ist, dass
die Sinusfunktion eine ganze Funktion ist, während die Tangensfunktion
Pole besitzt. Der Eisensteinsche Beweis des biquadratischen
Reziprozitätsgesetzes mit Hilfe elliptischer Funktionen ist daher
näher am Beweis mit der Tangensfunktion, denn elliptische Funktionen
haben zwangsläufig Pole.

Eisenstein spielt am Ende des \S~1 in \cite{Eis1} auf diesen Beweis
an; er schreibt nach der Bestimmung einer Konstanten $C$ in seinen
Kongruenzen:
\begin{quote}
  {\em Will man die Konstante $C$ vermeiden, muss man sich des Tangens
    anstatt des Sinus bedienen.}
\end{quote}
Koschmieder hat einen solchen Beweis mit Hilfe der Tangensfunktion
in  \cite{Ko2} gegeben und darüber am 20. September 1961 in
Oberwolfach vorgetragen; im Tagungsbericht steht dazu:
\begin{quote}
  {\em Eisenstein hat im Crelle Journal 29 (1845) einen Beweis des
    quadratischen Reziprozitätsgesetzes gegeben, bei dem er eine
    transzendente Hilfsfunktion benutzt, nämlich die Multiplikation
    des Sinus herangezogen hat. Er bemerkt dazu, daß der Beweis noch
    glatter liefe, wenn man statt des Sinus den Tangens gebrauchte. 
    Bisher scheint das nicht geschehen zu sein. Der Vortragende führt
    diesen Beweis.}
\end{quote}

\subsection*{Kommentare zum Beweis des kubischen Reziprozitätsgesetzes}
Den Beweis des kubischen Reziprozitätsgesetzes, den Eisenstein am
26. Januar 1845 gefunden hat, hat er ebenfalls nicht veröffentlicht; in
seinem Artikel \cite[\S~3]{Eis1} vom 13. Februar 1845 schreibt Eisenstein,
dass der Beweis dem des biquadratischen Reziprozitätsgesetzes ganz
analog sei. Eisensteins Beweisidee wurde später von zahlreichen
Autoren wieder aufgegriffen und ausgearbeitet, etwa von 
Hübler \cite{Huebler}, Dantscher\footnote{Die Namen Dantscher, Gegenbauer
  und Lewandowski zeigen das große Interesse, das den niederen
  Reziprozitätsgesetzen in Österreich entgegengebracht wurde. Auch Franz
  Mertens \cite{Mertens} hat sich mit analytischen Beweisen des kubischen und
  biquadratischen Reziprozitätsgesetzes beschäftigt, und Phänomene der
  Reziprozität spielen in der {\em Zahlentheorie} \cite{Aigner}
  Alexander Aigners eine große Rolle. Und dann ist
  da noch Emil Artin \ldots} \cite{Dantscher}, Gegenbauer \cite{Gegen},
Sbrana \cite{Sbrana}, Lewandowski \cite{Lewa}, Koschmieder \cite{Kosch},
Petr \cite{Petr}, Mel'nikov \cite{Mel} und Kubota \cite{Kubota}.

Koschmieder erwähnt in \cite{KoschDGL}, dass Lewandowski in
\cite{Lew9} auch das Reziprozitätsgesetz der 9ten Potenzreste mit
elliptischen Funktionen behandelt hat; ich habe aber weder diese Arbeit,
noch die darauf Bezug nehmende Arbeit von Georg Kantz \cite{Kantz29}
aus dem Jahre 1929 auffinden können.

\end{document}